\documentclass[noamsfonts]{amsproc}

\title[Finite Crystals and Paths]{Finite Crystals and Paths}

\author[G. Hatayama, Y. Koga, A. Kuniba, M. Okado and T. Takagi]{
Goro Hatayama, Yoshiyuki Koga, Atsuo Kuniba,\\
Masato Okado and Taichiro Takagi}

\address{
{\rm Goro Hatayama}\\
Institute of Physics\\
Graduate School of Science\\
University of Tokyo\\
Komaba, Tokyo 153-8902\\
Japan}
\address{
{\rm Yoshiyuki Koga}\\
Department of Mathematics\\
Graduate School of Science\\
Osaka University\\
Toyonaka, Osaka 560-0043\\
Japan}
\address{
{\rm Atsuo Kuniba}\\
Institute of Physics\\
Graduate School of Arts and Sciences\\
University of Tokyo\\
Komaba, Tokyo 153-8902\\
Japan}
\address{
{\rm Masato Okado}\\
Department of Informatics and Mathematical Science\\
Graduate School of Engineering Science\\
Osaka University\\
Toyonaka, Osaka 560-8531\\
Japan}
\address{
{\rm Taichiro Takagi}\\
Department of Mathematics and Physics\\
National Defense Academy\\
Yokosuka 239-8686\\
Japan}

\dedicatory{Dedicated to Professor Tetsuji Miwa on his fiftieth birthday}

\newtheorem{thm}{Donotwrite}[section]
\newtheorem{theorem}[thm]{Theorem}
\newtheorem{proposition}[thm]{Proposition}
\newtheorem{definition}[thm]{Definition}
\newtheorem{lemma}[thm]{Lemma}
\newtheorem{corollary}[thm]{Corollary}
\newtheorem{remark}[thm]{Remark}

\numberwithin{equation}{section}

\newfont{\germ}{eufm10}
\newfont{\slsmall}{cmsl8}
\newfont{\bfsl}{cmbxsl10}

\def\Bmin{B_{\min}}
\def\bb{\mbox{\bfsl b}}
\def\bt{\tilde{b}}
\def\C{{\mathcal C}}
\def\Cfin{\C^{fin}}
\def\Ch{\C^h}
\def\cd{\cdots}
\def\et#1{\tilde{e}_{#1}}
\def\ft#1{\tilde{f}_{#1}}
\def\geh{\goth{g}}
\def\goth#1{\mbox{\germ #1}}
\def\La{\Lambda}
\def\la{\lambda}
\def\lev{\mbox{\sl lev}\,}
\def\ol{\overline}
\def\ot{\otimes}
\def\P{{\mathcal P}}
\def\PB{{\mathcal P}(\pb,B)}
\def\pb{\mbox{\bfsl p}}
\def\Pcl{P_{cl}}
\def\Pcll{(P_{cl}^+)_l}
\def\Proof{\noindent{\sl Proof.}\quad}
\def\qed{~\rule{1mm}{2.5mm}}

\def\Uq{U_q(\geh)}
\def\Uqp{U'_q(\geh)}
\def\veps{\varepsilon}
\def\vphi{\varphi}
\def\wt{\mbox{\sl wt}\,}
\def\wts{\mbox{\slsmall wt}\,}
\def\Z{{\bf Z}}
\def\Zn{\Z_{\ge0}}

\begin{document}

\begin{abstract}
We consider a category of finite crystals of a quantum affine algebra
whose objects are not necessarily perfect, and set of paths, semi-infinite
tensor product of an object of this category with a certain boundary
condition. It is shown that the set of paths is isomorphic to a direct
sum of infinitely many, in general, crystals of integrable highest weight
modules. We present examples from $C_n^{(1)}$ and $A_{n-1}^{(1)}$, in which
the direct sum becomes a tensor product as suggested from the Bethe Ansatz.
\end{abstract}

\maketitle

\section{Introduction}

The main object of this note is to define a set of paths from a 
{\em finite} crystal $B$, which is not necessarily perfect, and
investigate its crystal structure. The set of paths $\PB$ is,
roughly speaking, a subset of the semi-infinite tensor product
$\cd\ot B\ot\cd\ot B\ot B$ with a certain boundary condition related
to $\pb$. If $B$ is perfect, it is known \cite{KMN1} that as crystals,
$\PB$ is isomorphic to the crystal base $B(\la)$ of an integrable highest 
weight module with highest weight $\la$ of the quantum affine algebra
$\Uq$. While trying to generalize this notion, we had two examples in mind:
(a) $\geh=C_n^{(1)},B=B^{1,l}\,(l:\mbox{odd})$; 
(b) $\geh=A_{n-1}^{(1)},B=B^{1,l}\ot B^{1,m}\,(l\ge m)$. 
For this parametrization
of finite crystals, we refer to \cite{HKOTY}. $B^{1,l}$ stands for the 
crystal base of an irreducible finite-dimensional $\Uqp$-module. 
In case (a) (resp. (b)) this finite-dimensional module is isomorphic 
to $V_{l\ol{\La}_1}\oplus V_{(l-2)\ol{\La}_1}\oplus\cd\oplus V_{\ol{\La}_1}$
(resp. $V_{l\ol{\La}_1}$) as $U_q(\ol{\geh})$-module, where $V_\la$ is 
the irreducible finite-dimensional module with highest weight $\la$.
In both cases $B$ is not perfect except when $l=m$ in (b).
For precise treatment see section \ref{subsec:ex-1} for (a) and
\ref{subsec:ex-2} for (b).

Let us consider case (a) first. When $l=1$ it has already been known
\cite{DJKMO} that the formal character of $\P(\pb,B^{1,1})$ for suitable
$\pb$ agrees with that of the irreducible highest weight 
$A_{2n-1}^{(1)}$-module with fundamental highest weight $\La_i$ regarded
as $C_n^{(1)}$-module via the natural embedding $C_n^{(1)}
\hookrightarrow A_{2n-1}^{(1)}$. On the other hand, the Bethe Ansatz 
suggests \cite{Ku} that $\P(\pb,B^{1,l})$ is equal to  $B(\la)\ot
\P(\pb^\dagger,B^{1,1})$ for suitable $\pb, \pb^\dagger$ and 
a level $\frac{l-1}{2}$ dominant integral weight $\la$ at the level of 
the Virasoro central charge.

Let us turn to case (b). In \cite{HKMW}
the $U'_q(\widehat{sl}_2)$-invariant integrable vertex model with 
alternating spins is considered. To translate the physical states and
operators of this model into the language of representation theory
of the quantum affine algebra $U_q(\widehat{sl}_2)$, they considered 
a set of paths with alternating spins and showed that it is isomorphic 
to the tensor product of crystals with highest weights.
Another appearance of example (b) can be found in \cite{HKKOTY}. They
considered the inductive limit of $(B^{1,l})^{\ot L_1}\ot
(B^{1,m})^{\ot L_2}$ when $L_1,L_2\rightarrow\infty,L_1\equiv r_1,
L_1+L_2\equiv r_2$ (mod $n$), and showed that there is a weight 
preserving bijection between the limit and $B((l-m)\La_{r_1})\ot 
B(m\La_{r_2})$. Since there is a natural isomorphism $B^{1,l}\ot B^{1,m}
\simeq B^{1,m}\ot B^{1,l}$, the above result claims that 
$\P(\pb,B^{1,l}\ot B^{1,m})$ for suitable $\pb$ is bijective to 
$B((l-m)\La_{r_1})\ot B(m\La_{r_2})$ with weight preserved.
These results are consistent with the earlier Bethe ansatz 
calculations on ``mixed spin" models 
\cite{AM,DMN}.

If we forget about the degree of the null root $\delta$ from weight,
this phenomenon is explained using the theory of crystals with
core \cite{KK}. (See also \cite{HKMW} section 3.2.)
Let $\{B_k\}_{k\ge1}$ be a coherent family of perfect crystals
and $B'_m$ be a perfect crystal of level $m$. 
Fix $l$ such that $l\ge m$ and take dominant integral weights $\la$ 
and $\mu$ of level $l-m$ and $m$. Then there exists an isomorphism of 
crystals:
\begin{eqnarray*}
B(\la)\ot B(\mu)&\simeq&B(\sigma\la)\ot B_{l-m}\ot B(\sigma'\mu)\ot B'_m\\
&\simeq& B(\sigma\la)\ot B(\sigma\sigma'\mu)\ot(B_l\ot B'_m),
\end{eqnarray*}
where $\sigma$ and $\sigma'$ are automorphisms on the weight lattice $P$
related to $\{B_k\}_{k\ge1}$ and $B'_m$. Iterating this isomorphism 
infinitely many times, we can expect
\[
\P(\pb^{(\la,\mu)},B_l\ot B'_m)\simeq B(\la)\ot B(\mu)
\]
as $P/\Z\delta$-weighted crystals with suitable $\pb^{(\la,\mu)}$.

In both cases (a),(b) we have illustrated above, what we expect is an
isomorphism of $P$-weighted crystals of the following type:
\begin{equation} \label{eq:intro}
\PB\simeq B(\la)\ot\P(\pb^\dagger,B^\dagger)
\end{equation}
and we shall prove it in this paper. First we examine the crystal
structure of $\PB$ and show it is isomorphic to a direct sum of
$B(\la)$'s. Therefore, the structure of $\PB$ is completely determined
by the set of highest weight elements. In the LHS of (\ref{eq:intro}),
such set $\PB_0$ is easy to describe,
and in the RHS, this set turns out to be the set of restricted paths
$\P^{(\la)}(\pb^\dagger,B^\dagger)$, which is familiar to the people 
in solvable lattice models. Thus establishing a weight preserving 
bijection between $\PB_0$ and $\P^{(\la)}(\pb^\dagger,B^\dagger)$
directly, we can show (\ref{eq:intro}).

\section{Crystals}

\subsection{Notation}
Let $\geh$ be an affine Lie algebra. We denote by $I$ the index set of its
Dynkin diagram. Note that $0$ is included in $I$. Let $\alpha_i,h_i,\La_i$
($i\in I$) be the simple roots, simple coroots, fundamental weights for 
$\geh$. Let $\delta=\sum_{i\in I}a_i\alpha_i$ denote the standard null 
root, and $c=\sum_{i\in I}a_i^\vee h_i$ the canonical central element,
where $a_i,a_i^\vee$ are positive integers as in \cite{Kac}. We assume 
$a_0=1$. Let $P=\bigoplus_{i\in I}\Z\La_i\oplus\Z\delta$ be the weight 
lattice, and set $P^+=\sum_{i\in I}\Zn\La_i\oplus\Z\delta$.

Let $\Uq$ be the quantum affine algebra associated to $\geh$. For the
definition of $\Uq$ and its Hopf algebra structure, see e.g. section
2.1 of \cite{KMN1}. For $J\subset I$ we denote by $U_q(\geh_J)$ 
the subalgebra of $\Uq$ generated by $e_i,f_i,t_i$ ($i\in J$).
In particular, $U_q(\geh_{I\setminus\{0\}})$ is identified with the quantized
enveloping algebra for the simple Lie algebra whose Dynkin diagram is 
obtained by deleting the 0 vertex from that of $\geh$.
We also consider the quantum affine algebra without
derivation $\Uqp$. As its weight lattice, the classical weight lattice
$\Pcl=P/\Z\delta$ is needed. We canonically identify $\Pcl$ with 
$\bigoplus_{i\in I}\Z\La_i\subset P$. For the precise treatment,
see section 3.1 of \cite{KMN1}.
We further define the following subsets of 
$\Pcl$: $\Pcl^0=\{\la\in \Pcl\mid \langle\la,c\rangle=0\}$,
$\Pcl^+=\{\la\in \Pcl\mid \langle\la,h_i\rangle\ge0\mbox{ for any }i\}$,
$\Pcll=\{\la\in \Pcl^+\mid \langle\la,c\rangle=l\}$. 
For $\la,\mu\in\Pcl$, we write $\la\ge\mu$ to mean $\la-\mu\in\Pcl^+$.

\subsection{Crystals and crystal bases}
We summarize necessary facts in crystal theory. Our basic references are 
\cite{K1}, \cite{KMN1} and \cite{AK}. 

A crystal $B$ is a set $B$ with the maps
\[
\et{i},\ft{i}: B\sqcup\{0\}\longrightarrow B\sqcup\{0\}
\]
satisfying the following properties:
\begin{itemize}
\item[] $\et{i}0=\ft{i}0=0$,
\item[] for any $b$ and $i$, there exists $n>0$ such that 
$\et{i}^nb=\ft{i}^nb=0$,
\item[] for $b,b'\in B$ and $i\in I$, $\ft{i}b=b'$ if and only if 
$b=\et{i}b'$.
\end{itemize}
If we want to emphasize $I$, $B$ is called an $I$-crystal.
A crystal can be regarded as a colored oriented graph by defining
\[
b\stackrel{i}{\longrightarrow}b'\quad\Longleftrightarrow\quad \ft{i}b=b'.
\]
For an element $b$ of $B$ we set
\[
\veps_i(b)=\max\{n\in\Zn\mid\et{i}^nb\neq0\},\quad
\vphi_i(b)=\max\{n\in\Zn\mid\ft{i}^nb\neq0\}.
\]
We also define a $P$-weighted crystal. It is a crystal with the weight
decomposition $B=\sqcup_{\la\in P}B_\la$ such that 
\begin{eqnarray}
&&\et{i}B_\la\subset B_{\la+\alpha_i}\sqcup\{0\},\quad
\ft{i}B_\la\subset B_{\la-\alpha_i}\sqcup\{0\},\\
&&\langle h_i,\wt b\rangle=\vphi_i(b)-\veps_i(b). \label{eq:i-wt}
\end{eqnarray}
Set
\[
\veps(b)=\sum_{i\in I}\veps_i(b)\La_i,\quad
\vphi(b)=\sum_{i\in I}\vphi_i(b)\La_i.
\]
Then (\ref{eq:i-wt}) is equivalent to $\vphi(b)-\veps(b)=\wt b$. 
$\Pcl$-weighted crystal is defined similarly.

For two weighted crystals $B_1$ and $B_2$, the tensor product $B_1\ot B_2$
is defined.
\[
B_1\ot B_2=\{b_1\ot b_2\mid b_1\in B_1,b_2\in B_2\}.
\]
The actions of $\et{i}$ and $\ft{i}$ are defined by
\begin{eqnarray}
\et{i}(b_1\ot b_2)&=&\left\{
\begin{array}{ll}
\et{i}b_1\ot b_2&\mbox{ if }\vphi_i(b_1)\ge\veps_i(b_2)\\
b_1\ot \et{i}b_2&\mbox{ if }\vphi_i(b_1) < \veps_i(b_2),
\end{array}\right. \label{eq:ot-e}\\
\ft{i}(b_1\ot b_2)&=&\left\{
\begin{array}{ll}
\ft{i}b_1\ot b_2&\mbox{ if }\vphi_i(b_1) > \veps_i(b_2)\\
b_1\ot \ft{i}b_2&\mbox{ if }\vphi_i(b_1)\le\veps_i(b_2).
\end{array}\right. \label{eq:ot-f}
\end{eqnarray}
Here $0\ot b$ and $b\ot0$ are understood to be $0$.
$\veps_i,\vphi_i$ and $\wt$ are given by
\begin{eqnarray}
\veps_i(b_1\ot b_2)&=&
\max(\veps_i(b_1),\veps_i(b_1)+\veps_i(b_2)-\vphi_i(b_1)),\label{eq:ot-eps}\\
\vphi_i(b_1\ot b_2)&=&
\max(\vphi_i(b_2),\vphi_i(b_1)+\vphi_i(b_2)-\veps_i(b_2)),\label{eq:ot-phi}\\
\wt(b_1\ot b_2)&=&\wt b_1+\wt b_2.
\end{eqnarray}

\begin{definition}[\cite{AK}]
We say a $P$ (or $\Pcl$)-weighted crystal is regular, if for any $i,j\in I$
($i\ne j$), $B$ regarded as $\{i,j\}$-crystal is a disjoint union of 
crystals of integrable highest weight modules over $U_q(\geh_{\{i,j\}})$.
\end{definition}

Crystal is a notion obtained by abstracting the properties of crystal 
bases \cite{K1}. Let $V(\la)$ be the integrable highest weight $\Uq$-module
with highest weight $\la\in P^+$ and highest weight vector $u_\la$. It is
shown in \cite{K1} that $V(\la)$ has a crystal base $(L(\la),B(\la))$. 
We regard $u_\la$ as an element of $B(\la)$ as well. $B(\la)$ is a regular
$P$-weighted crystal. A finite-dimensional integrable $\Uqp$-module $V$
does not necessarily have a crystal base. If $V$ has a crystal base $(L,B)$,
then $B$ is a regular $\Pcl^0$-weighted crystal with finitely many 
elements.

Let $W$ be the affine Weyl group associated to $\geh$, and $s_i$ be the 
simple reflection corresponding to $\alpha_i$. $W$ acts on any regular 
crystal $B$ \cite{K2}. The action is given by
\[
S_{s_i}b=\left\{
\begin{array}{ll}
\ft{i}^{\langle h_i,\wts b\rangle}b\quad
&\mbox{ if }\langle h_i,\wt b\rangle\ge0\\
\et{i}^{-\langle h_i,\wts b\rangle}b\quad
&\mbox{ if }\langle h_i,\wt b\rangle\le0.
\end{array}
\right.
\]
An element $b$ of $B$ is called $i$-{\em extremal} if $\et{i}b=0$ or 
$\ft{i}b=0$. $b$ is called {\em extremal} if $S_wb$ is $i$-extremal
for any $w\in W$ and $i\in I$. 

\begin{definition}[\cite{AK} Definition 1.7] \label{def:simple}
Let $B$ be a regular $\Pcl^0$-weighted crystal with finitely many 
elements.
We say $B$ is simple if it satisfies
\begin{itemize}
\item[(1)] There exists $\la\in\Pcl^0$ such that the weights of $B$ are
in the convex hull of $W\la$.
\item[(2)] $\sharp B_\la=1$.
\item[(3)] The weight of any extremal element is in $W\la$.
\end{itemize}
\end{definition}

\begin{remark}
Let $B$ be a regular $\Pcl^0$-weighted crystal with finitely many elements. 
We have the following criterion for simplicity. Let $B(\la)$ denote 
the crystal base of the irreducible highest weight 
$U_q(\geh_{I\setminus\{0\}})$-module with
highest weight $\la$. If $B$ decomposes into 
$B\simeq\bigoplus_{j=0}^m B(\la_j)$
as $U_q(\geh_{I\setminus\{0\}})$-crystal and $\la_j$ satisfies 
\begin{itemize}
\item[(1)] $\la_j\in\la_0+\sum_{i\ne0}\Z_{\le0}\alpha_i$ and $\la_j\ne\la_0$
for any $j\ne0$,
\item[(2)] The highest weight element of $B(\la_j)$ is not $0$-extremal 
for any $j\ne0$,
\end{itemize}
then $B$ is simple.
\end{remark}

\begin{proposition}[\cite{AK} Lemma 1.9 \& 1.10] \label{prop:simple}
Simple crystals have the following properties.
\begin{itemize}
\item[(1)] A simple crystal is connected.
\item[(2)] The tensor product of simple crystals is also simple.
\end{itemize}
\end{proposition}

\subsection{Category $\Cfin$}
Let $B$ be a regular $\Pcl^0$-weighted crystal with finitely many 
elements. For $B$ we introduce the {\em level} of $B$ by 
\[
\lev B=\min\{\langle c,\veps(b)\rangle\mid b\in B\}\in\Zn.
\]
Note that $\langle c,\veps(b)\rangle=\langle c,\vphi(b)\rangle$ for any 
$b\in B$. We also set $\Bmin=\{b\in B\mid\langle c,\veps(b)\rangle=\lev B\}$
and call an element of $\Bmin$ {\em minimal}.

\begin{definition} \label{def:C-fin}
We denote by $\Cfin(\geh)$ (or simply $\Cfin$) the category of crystal
$B$ satisfying the following conditions:
\begin{itemize}
\item[(1)] $B$ is a crystal base of a finite-dimensional $\Uqp$-module.
\item[(2)] $B$ is simple.
\item[(3)] For any $\la\in\Pcl^+$ such that $\langle c,\la\rangle\ge
           \lev B$, there exists $b\in B$ satisfying $\veps(b)\le\la$.
           It is also true for $\vphi$.
\end{itemize}
\end{definition}

We call an object of $\Cfin(\geh)$ {\em finite crystal}.

\begin{remark}
\begin{itemize}
\item[(i)] Condition (1) implies $B$ is a regular $\Pcl^0$-weighted
           crystal with finitely many elements.
\item[(ii)] Set $l=\lev B$. Condition (3) implies that the maps $\veps$
           and $\vphi$ from $\Bmin$ to $\Pcll$ are surjective.
           (cf. (4.6.5) in \cite{KMN1})
\item[(iii)] Practically, one has to check condition (3) only for 
           $\la\in\Pcl^+$ such that there is no $i\in I$ satisfying
           $\la-\La_i\ge0$ and $\langle c,\la-\La_i\rangle\ge\lev B$.
           In particular, if $a^\vee_i=1$ for any $i\in I$
           ($\geh=A_n^{(1)},C_n^{(1)}$), the surjectivity of $\veps$ and
           $\vphi$ assures (3).
\item[(iv)] The authors do not know a crystal satisfying (1) and (2), but
           not satisfying (3).
\end{itemize}
\end{remark}

Let $B_1$ and $B_2$ be two finite crystals. Definition
\ref{def:C-fin} (1) and the existence of the universal $R$-matrix assures
that we have a natural isomorphism of crystals.
\begin{equation} \label{eq:iso}
B_1\ot B_2\simeq B_2\ot B_1.
\end{equation}

The following lemma is immediate.

\begin{lemma} \label{lem:tensor}
Let $B_1,B_2$ be finite crystals.
\begin{itemize}
\item[(1)] $\lev(B_1\ot B_2)=\max(\lev B_1,\lev B_2)$.
\item[(2)] If $\lev B_1\ge\lev B_2$, then $(B_1\ot B_2)_{\min}=
           \{b_1\ot b_2\mid b_1\in(B_1)_{\min},$
           $\vphi_i(b_1)\ge\veps_i(b_2)\mbox{ for any }i\}$.
\item[(3)] If $\lev B_1\le\lev B_2$, then $(B_1\ot B_2)_{\min}=
           \{b_1\ot b_2\mid b_2\in(B_2)_{\min},$
           $\vphi_i(b_1)\le\veps_i(b_2)\mbox{ for any }i\}$.
\end{itemize}
\end{lemma}

$\Cfin(\geh)$ forms a tensor category.

\begin{proposition}
If $B_1$ and $B_2$ are objects of $\Cfin(\geh)$, then $B_1\ot B_2$ is
also an object of $\Cfin(\geh)$.
\end{proposition}

\Proof
We need to check the conditions in Definition \ref{def:C-fin} for 
$B_1\ot B_2$. (1) is obvious and (2) follows from Proposition 
\ref{prop:simple} (2).

Let us prove condition (3) for $\veps$. Set $l_1=\lev B_1,l_2=\lev B_2$.
Using (\ref{eq:iso}) if necessary, we can assume $l_1\ge l_2$.
Thus we have $\lev B_1\ot B_2=l_1$. For any $\la\in\Pcl^+$ such that
$\langle c,\la\rangle\ge l_1$, one can take $b_1\in B_1$ satisfying
$\veps(b_1)\le\la$. Since $\langle c,\vphi(b_1)\rangle\ge l_1\ge l_2$,
one can take $b_2\in B_2$ satisfying $\veps(b_2)\le \vphi(b_1)$.
In view of (\ref{eq:ot-eps}) one has $\veps(b_1\ot b_2)=\veps(b_1)\le\la$.

For the proof of $\vphi$, repeat a similar exercise for $B_2\ot B_1
(\simeq B_1\ot B_2)$ using (\ref{eq:ot-phi}).
\qed

\subsection{Category $\Ch$} \label{subsec:C-h}
If an element $b$ of a crystal $B$ satisfies $\et{i}b=0$ for any $i$, we 
call it a {\em highest weight} element.

\begin{definition}
We denote by $\Ch(I,P)$ (or simply $\Ch$) the category of regular $P$-weighted
crystal $B$ satisfying the following condition:
\begin{itemize}
\item[] For any $b\in B$, there exist $l\ge0,i_1,\cd,i_l\in I$ such that 
$b'=\et{i_1}\cd\et{i_l}b\in B$ is a highest weight element.
\end{itemize}
\end{definition}
Clearly, $\Ch(I,P)$ forms a tensor category. 

\begin{proposition}[\cite{KMN1} Proposition 2.4.4]
An object of $\Ch(I,P)$ is isomorphic to a direct sum (disjoint union) 
of crystals $B(\la)$ ($\la\in P^+$) of integrable highest weight $\Uq$-modules.
\end{proposition}

Let $O$ be an object of $\Ch(I,P)$. By $O_0$ we mean the set of highest 
weight elements in $O$. Suppose that $O_0=\{b_j\mid j\in J\}$ and 
$\wt b_j=\la_j\in P^+$, then from the above proposition we have an
isomorphism
\[
O\simeq\bigoplus_{j\in J}B(\la_j)\quad\mbox{as $P$-weighted crystals}.
\]
$J$ can be an infinite set.

The following lemma is standard.

\begin{lemma} \label{lem:standard}
Let $B_1,B_2$ be weighted crystals. Then $b_1\ot b_2\in B_1\ot B_2$ is a 
highest weight element, if and only if $b_1$ is a highest weight element
and $\et{i}^{\langle h_i,\wts b_1\rangle+1}b_2=0$ for any $i$.
\end{lemma}

Let $O$ be an object of $\Ch(I,P)$. From this lemma we have the following
bijection.
\begin{eqnarray*}
(B(\la)\ot O)_0&\quad\longrightarrow\quad&
O^{\le\la}:=\{b\in O\mid\et{i}^{\langle h_i,\la\rangle+1}b=0
\mbox{ for any }i\}\\
u_\la\ot b&\mapsto&b.
\end{eqnarray*}
Note that $O^{\le0}=O_0$.

\section{Paths}

In this section we construct a set of paths from a finite crystal
and consider its structure.

\subsection{Energy function}
Let us recall the energy function used in \cite{NY} to identify the 
Kostka-Foulkes polynomial with a generating function over classically 
restricted paths.

Let $B_1$ and $B_2$ be two finite crystals. 
Suppose $b_1\ot b_2\in B_1\ot B_2$ is mapped to $\bt_2\ot\bt_1
\in B_2\ot B_1$ under the isomorphism (\ref{eq:iso}). A $\Z$-valued function
$H$ on $B_1\ot B_2$ is called an {\em energy function} if for any $i$
and $b_1\ot b_2\in B_1\ot B_2$ such that $\et{i}(b_1\ot b_2)\neq0$,
 it satisfies
\begin{eqnarray}
H(\et{i}(b_1\ot b_2))&=H(b_1\ot b_2)+1
&\mbox{ if }i=0,\vphi_0(b_1)\geq\veps_0(b_2),
\nonumber\\
&&\phantom{\mbox{ if }}\vphi_0(\bt_2)\geq\veps_0(\bt_1),
\nonumber\\
&=H(b_1\ot b_2)-1
&\mbox{ if }i=0,\vphi_0(b_1)<\veps_0(b_2),
\nonumber\\
&&\phantom{\mbox{ if }}\vphi_0(\bt_2)<\veps_0(\bt_1),
\nonumber\\
&\hspace{-6mm}=H(b_1\ot b_2)&\mbox{ otherwise}.\label{eq:e-func}
\end{eqnarray}
When we want to emphasize $B_1\ot B_2$, we write $H_{B_1B_2}$ for $H$.
The existence of such function can be shown in a similar manner to 
section 4 of \cite{KMN1} based on the existence of {\em combinatorial 
$R$-matrix}. The energy function is unique up to additive constant, since
$B_1\ot B_2$ is connected. By definition, $H_{B_1B_2}(b_1\ot b_2)
=H_{B_2B_1}(\bt_2\ot\bt_1)$.

If the tensor product $B_1\ot B_2$ is homogeneous, i.e., $B_1=B_2$,
we have $\bt_2=b_1,\bt_1=b_2$. Thus (\ref{eq:e-func}) is rewritten as
\begin{eqnarray}
H(\et{i}(b_1\ot b_2))&=H(b_1\ot b_2)+1
&\mbox{ if }i=0,\vphi_0(b_1)\geq\veps_0(b_2),\nonumber\\
&=H(b_1\ot b_2)-1
&\mbox{ if }i=0,\vphi_0(b_1)<\veps_0(b_2),\nonumber\\
&\hspace{-6mm}=H(b_1\ot b_2)&\mbox{ if }i\ne0.\label{eq:e-func*}
\end{eqnarray}

The following proposition, which is shown by case-by-case checking,
reduces the energy function of a tensor product 
to that of each component.

\begin{proposition} \label{prop:HofTensor}
Set $B=B_1\ot B_2$, then
\begin{eqnarray*}
H_{BB}((b_1\ot b_2)\ot(b'_1\ot b'_2))
&=&H_{B_1B_2}(b_1\ot b_2)+H_{B_1B_1}(\bt_1\ot b'_1)\\
&&+H_{B_2B_2}(b_2\ot \bt'_2)+H_{B_1B_2}(b'_1\ot b'_2).
\end{eqnarray*}
Here $\bt_1,\bt'_2$ are defined as
\begin{eqnarray*}
B_1\ot B_2&\simeq&B_2\ot B_1\\
b_1\ot b_2&\mapsto&\bt_2\ot\bt_1\\
b'_1\ot b'_2&\mapsto&\bt'_2\ot\bt'_1.
\end{eqnarray*}
\end{proposition}

\begin{remark}
Decomposition of the energy function is not unique.
For instance, the following also gives such decomposition.
\begin{eqnarray*}
H_{BB}((b_1\ot b_2)\ot(b'_1\ot b'_2))
&=&H_{B_2B_1}(b_2\ot b'_1)+H_{B_1B_1}(b_1\ot\check{b}'_1)\\
&&+H_{B_2B_2}(\check{b}_2\ot b'_2)+H_{B_1B_2}(\check{b}'_1\ot\check{b}_2),
\end{eqnarray*}
where 
\begin{eqnarray*}
B_2\ot B_1&\simeq&B_1\ot B_2\\
b_2\ot b'_1&\mapsto&\check{b}'_1\ot\check{b}_2.
\end{eqnarray*}
\end{remark}

\subsection{Set of paths $\PB$}
We shall define a set of paths from any finite crystal in $\Cfin$ 
imitating the construction in section 4 of \cite{KMN1} from a perfect crystal.

\begin{definition}
An element 
$\pb=\cd\ot\bb_j\ot\cd\ot\bb_2\ot\bb_1$ of the semi-infinite tensor product
of $B$ is called a reference path if it satisfies $\bb_j\in \Bmin$ and
$\vphi(\bb_{j+1})=\veps(\bb_j)$ for any $j\ge1$.
\end{definition}

\begin{definition}
Fix a reference path $\pb=\cd\ot\bb_j\ot\cd\ot\bb_2\ot\bb_1$. We define 
a set of paths $\PB$ by
\[
\PB=\{p=\cd\ot b_j\ot\cd\ot b_2\ot b_1\mid b_j\in B,b_k=\bb_k\mbox{ for }
k\gg1\}.
\]
\end{definition}
An element of $\PB$ is called a {\em path}. For convenience we denote
$b_k$ by $p(k)$ and $\cd\ot b_{k+2}\ot b_{k+1}$ by $p[k]$ for 
$p=\cd\ot b_j\ot\cd\ot b_2\ot b_1$.

\begin{definition} \label{def:E-and-weight}
For a path $p\in\PB$, set
\begin{eqnarray*}
E(p)&=&\sum_{j=1}^\infty j(H(p(j+1)\ot p(j))-H(\pb(j+1)\ot\pb(j))),\\
W(p)&=&\vphi(\pb(1))+\sum_{j=1}^\infty (\wt p(j)-\wt \pb(j))
-E(p)\delta.
\end{eqnarray*}
$E(p)$ and $W(p)$ are called the energy and weight of $p$.
\end{definition}

We distinguish $W(p)\in P$ from 
$\wt p=\vphi(\pb(1))+\sum_{j=1}^\infty (\wt p(j)-\wt \pb(j))\in\Pcl$.

\begin{remark}
\begin{itemize}
\item[(i)] If $B$ is perfect, the set of reference paths is bijective
           to $\Pcll$, where $l=\lev B$. For $\la\in\Pcll$ take a unique
           $\bb_1\in\Bmin$ such that $\vphi(\bb_1)=\la$. The condition
           $\vphi(\bb_{j+1})=\veps(\bb_j)$ fixes 
           $\pb=\cd\ot\bb_j\ot\cd\ot\bb_1$ uniquely.
\item[(ii)] In \cite{KMN1} $\pb$ is called a ground state path,
           since $E(p)\ge E(\pb)$ for any $p\in\PB$. But if $B$ is not
           perfect, it is no longer true in general.
\end{itemize}
\end{remark}

The following theorem is essential for our consideration below.

\begin{theorem} \label{th:C-h}
Assume $\mbox{\sl rank }\geh>2$. Then $\PB$ is an object of $\Ch$.
\end{theorem}

\Proof
Assume $\et{i}p=\cd\ot\et{i}b_j\ot\cd\ot b_1\ne0$. Note that 
$E(\et{i}p)=E(p)-\delta_{i0}$ and $\wt\et{i}b_j=\wt b_j+\alpha_i-
\delta_{i0}\delta\in\Pcl$. By Definition \ref{def:E-and-weight}
it is immediate to see $\PB$ is a $P$-weighted crystal. Thus one has to
check the following:
\begin{itemize}
\item[(i)] If for any $i,j\in I$ ($i\ne j$), $\PB$ regarded as 
           $\{i,j\}$-crystal is a disjoint union of crystals of integrable
           highest weight modules over $U_q(\geh_{\{i,j\}})$.
\item[(ii)] For any $p\in\PB$, there exist $l\ge0,i_1,\cd,i_l\in I$
           such that $p'=\et{i_1}\cd\et{i_l}p\in\PB$ is a highest weight
           element.
\end{itemize}

We prove (i) first. For $p\in\PB$ take $m,m'$ such that $p(k)=\pb(k)$ 
for $k>m$ and $m'\gg m$. Note that if $\ft{i_N}\cd\ft{i_1}p[m]=
p[m']\ot b'_{m'}\ot\cd\ot b'_{m+1}$, then $b'_k=\pb(k)$ for $k>m+N$. From the 
assumption, $U_q(\geh_{\{i,j\}})$ is the quantized enveloping algebra 
associated to a finite-dimensional Lie algebra. Since $B$ is regular,
the connected component containing $p[m]$, as $\{i,j\}$-crystal, can be
considered to be in $B(\vphi(p[m']))\ot B^{\ot(m'-m)}$. 
Since $\veps(p[m])=0$, we can regard $p[m]$ as highest 
weight element of some $\{i,j\}$-crystal $B_0$ which is isomorphic
to the crystal of an integrable highest weight $U_q(\geh_{\{i,j\}})$-module.
Hence $p$ is contained in a component of the $\{i,j\}$-crystal
$B_0\ot B^{\ot m}$, which is a disjoint union of crystals of integrable 
highest weight $U_q(\geh_{\{i,j\}})$-modules.

To prove (ii) for $p=\cd\ot b_k\ot\cd\ot b_1\in\PB$, we take the minimum
integer $m$ such that $p'=p[m]$ is a highest weight element. We prove by
induction on $m$.

First let us show that there exist $l\ge0,i_1,\cd,i_l\in I$ such that
$\et{i_1}\cd\et{i_l}(p'\ot b_m)$ is a highest weight element. 
The proof is essentially the same as a part of that of Theorem 4.4.1
in \cite{KMN1}. Nevertheless we repeat it for the sake of self-containedness.
Suppose that there does not exist such $i_1,\cd,i_l$. Then there exists an 
infinite sequence $\{i_\nu\}$ in $I$ such that 
\[
\et{i_k}\cd\et{i_1}(p'\ot b_m)\neq0.
\]
Since $\et{i_k}\cd\et{i_1}(p'\ot b_m)=p'\ot\et{i_k}\cd\et{i_1}b_m$ and
$B$ is a finite set, there exists $b^{(1)}\in B$ and $j_1,\cd,j_l$ such that
\[
p'\ot b^{(1)}=\et{j_l}\cd\et{j_1}(p'\ot b^{(1)}).
\]
Hence setting $b^{(\nu+1)}=\et{j_\nu}b^{(\nu)}$, we have
\[
\et{j_\nu}(p'\ot b^{(\nu)})=p'\ot b^{(\nu+1)}\mbox{ and }b^{(l+1)}=b^{(1)}.
\]
In view of (\ref{eq:ot-phi}) we have $\vphi_i(p')\ge\vphi_i(b_{m+1})$ for
any $i$. Thus by (\ref{eq:ot-e}) we have $\veps_{j_\nu}(b^{(\nu)})>
\vphi_{j_\nu}(p')\ge\vphi_{j_\nu}(b')$ for some $b'\in B$. Hence we have
\[
\et{j_\nu}(b'\ot b^{(\nu)})=b'\ot b^{(\nu+1)}.
\]
Therefore, from (\ref{eq:e-func*}), we have 
\[
H(b'\ot b^{(\nu+1)})=H(b'\ot b^{(\nu)})-\delta_{i_\nu0}.
\]
Hence $H(b'\ot b^{(l+1)})=H(b'\ot b^{(1)})-\sharp\{\nu\mid j_\nu=0\}$,
which implies there is no $\nu$ such that $j_\nu=0$. On the other hand,
$\sum_\nu\alpha_{j_\nu}=0\mbox{ mod }\Z\delta$ and hence $\sum_\nu
\alpha_{j_\nu}$ is a positive multiple of $\delta$, which contradicts 
$0\notin\{j_1,\cd,j_l\}$.

Now set $p''=p'\ot b_m(=p[m-1]),b''=b_{m-1}\ot\cd\ot b_1$. Notice that
for any $i\in I$ satisfying $\et{i}p''\neq0$, there exists $k\ge1$ such that
\[
\et{i}^k(p''\ot b'')=\et{i}p''\ot\et{i}^{k-1}b''.
\]
Therefore there exist $l\ge0,(i_1,k_1),\cd,(i_l,k_l)\in I\times\Z_{>0}$
such that 
\[
\et{i_1}^{k_1}\cd\et{i_l}^{k_l}p=\et{i_1}\cd\et{i_l}p''\ot
\et{i_1}^{k_1-1}\cd\et{i_l}^{k_l-1}b''
\]
and $\et{i_1}\cd\et{i_l}p''$ is a highest weight element. Now we can use
the induction assumption and complete the proof.
\qed

\begin{remark}
As seen in the proof, the theorem does not require the condition 
$\bb_j\in\Bmin$ for the reference path $\pb=\cd\ot\bb_j\ot\cd\ot\bb_1$.
\end{remark}

The following proposition describes the set of highest weight elements in
$\PB$.

\begin{proposition} \label{prop:P_0}
\[
\PB_0=\{p\in\PB\mid p(j)\in\Bmin,
\vphi(p(j+1))=\veps(p(j))\mbox{ for }\forall j\}.
\]
\end{proposition}

\Proof
Assume $p=\cd\ot b_j\ot\cd\ot b_1$ is a highest weight element. 
We prove the following by induction on $m$ in decreasing order.
\begin{itemize}
\item[(i)] $b_m\in\Bmin,\vphi(b_{m+1})=\veps(b_m)$
\item[(ii)] $\vphi(p[m-1])=\vphi(b_m)$
\end{itemize}
These conditions are satisfied for sufficiently large $m$.
{}From (ii) for $m+1$ we have $\vphi(p[m])=\vphi(b_{m+1})$.
{}From Lemma \ref{lem:standard} we see that $p[m]$ is a highest 
weight element and $\veps(b_m)\le\wt p[m]=\vphi(p[m])=\vphi(b_{m+1})$.
Combining this with (i) for $m+1$, we can conclude (i) for $m$.
For (ii) use (\ref{eq:ot-phi}).
\qed

As seen in the proof, we obtain

\begin{corollary} \label{cor:useful}
If $p\in\PB_0$, then $\wt p[j]=\vphi(p(j+1))$.
\end{corollary}

\subsection{Restricted paths}
When $B$ is perfect the set of {\em restricted} paths was defined in
\cite{DJO} and shown to be bijective to $(B(\la)\ot B(\mu))_0$ for some
$\la,\mu\in\Pcl^+$. Here we shall consider restricted paths for any finite
crystal $B$.

For $\la\in\Pcl^+$ and $p\in\PB$, we introduce a
sequence of weights $\{\la_j(p)\}_{j\ge0}$ by
\begin{eqnarray*}
&&\la_j(p)=\la+\vphi(p(j+1))\mbox{ for }j\gg1,\\
&&\la_{j-1}(p)=\la_j(p)+\wt p(j).
\end{eqnarray*}
Notice that this definition is well-defined by virtue of the property
of the reference path. In fact, $\la_j(p)=\la+\wt p[j]$.

\begin{definition}
For $\la\in\Pcl^+$ we define a subset $\P^{(\la)}(\pb,B)$ of $\PB$ by 
\[
\P^{(\la)}(\pb,B)=\{p\in\PB\mid
\et{i}^{\langle h_i,\la_j(p)\rangle+1}p(j)=0\mbox{ for }\forall i,j\}.
\]
\end{definition}

An element of $\P^{(\la)}(\pb,B)$ is called a {\em restricted} path.

\begin{proposition} \label{prop:restricted}
For $\la\in\Pcl^+$ we have 
\[
\PB^{\le\la}=\P^{(\la)}(\pb,B).
\]
\end{proposition}

\Proof
Assume $p=\cd\ot b_j\ot\cd\ot b_1\in\PB^{\le\la}$, which is equivalent to
saying $u_\la\ot p$ is a highest weight element. So is $u_\la\ot p[j]\ot b_j$
by Lemma \ref{lem:standard}. Using this lemma again we get 
$\veps(b_j)\le\wt(u_\la\ot p[j])=\la_j(p)$.

To show the inverse inclusion, assume $p=\cd\ot b_j\ot\cd\ot b_1\in
\P^{(\la)}(\pb,B)$. We prove $\veps(p[j])\le\la$ by induction on $j$ in
decreasing order. We know $\veps(p[j])=0$ for sufficiently large $j$.
Supposing $\veps(p[j])\le\la$ we immediately obtain $\veps(p[j]\ot b_j)
\le\la$ from (\ref{eq:ot-eps}) and the condition $\veps(b_j)\le\la_j(p)$.
\qed

As seen in the proof we have $\la_j(p)\in\Pcl^+$ and its level is
$\langle c,\la\rangle+\lev B$.

Combining the results in section \ref{subsec:C-h}, Theorem \ref{th:C-h}
and Proposition \ref{prop:restricted}, we obtain

\begin{theorem} \label{th:main}
Let $\PB$ and $\P(\pb^\dagger,B^\dagger)$ be two sets of paths.
If for certain $\la\in\Pcl^+$, there exists a bijection
\begin{eqnarray} 
\PB_0&\longrightarrow&\P^{(\la)}(\pb^\dagger,B^\dagger) \label{eq:P_0=P^la}\\
p&\mapsto&p^\dagger \nonumber
\end{eqnarray}
such that $W(p)=\la+W(p^\dagger)$, then we have an isomorphism of
$P$-weighted crystals
\[
\PB\simeq B(\la)\ot\P(\pb^\dagger,B^\dagger).
\]
They are isomorphic to a direct sum of crystals of integrable highest 
weight $\Uq$-modules, and their highest weight elements are parametrized by
(\ref{eq:P_0=P^la}).
\end{theorem}

\section{Examples}

We shall give two examples to which we can apply Theorem \ref{th:main}
efficiently.

\subsection{Example 1} \label{subsec:ex-1}
We present a useful proposition first.
Similar to $O^{\le\la}$ we define $B^{\le\la}$ for a finite crystal $B$
and $\la\in\Pcl^+$ by
\[
B^{\le\la}=\{b\in B\mid \et{i}^{\langle h_i,\la\rangle+1}b=0
\mbox{ for any }i\}.
\]
Note that if $\lev B=l$, then $\Bmin=\bigsqcup_{\la\in\Pcll}B^{\le\la}$.

\begin{proposition} \label{prop:case-1}
Let $B$ and $B^\dagger$ be finite crystals such that 
$\lev B\ge\lev B^\dagger$, and 
$\pb=\cd\ot\bb_j\ot\cd\ot\bb_1$ be a reference path for $B$. 
Suppose there exists a map $t:\Bmin\rightarrow B^\dagger$ satisfying
the following conditions:
\begin{itemize}
\item[(1)] For any $\mu\in\Pcll$ ($l=\lev B$), $t|_{B^{\le\mu}}$ is a 
           bijection onto $(B^\dagger)^{\le\mu}$.
\item[(2)] $\wt t(b)=\wt b$ for any $b\in\Bmin$.
\item[(3)] $H_{B^\dagger B^\dagger}(t(b_1)\ot t(b_2))=H_{BB}(b_1\ot b_2)$
           up to global additive constant for any $(b_1,b_2)\in\Bmin^2$ 
           such that $\vphi(b_1)=\veps(b_2)$.
\item[(4)] $\pb^\dagger=\cd\ot t(\bb_j)\ot\cd\ot t(\bb_1)$ is a reference
           path for $B^\dagger$.
\end{itemize}
Then setting $\la=\vphi(\bb_1)-\vphi(t(\bb_1))$, we have
\[
\PB\simeq B(\la)\ot\P(\pb^\dagger,B^\dagger).
\]
\end{proposition}

\Proof
Consider the following map.
\begin{eqnarray*}
\PB_0&\longrightarrow&\P(\pb^\dagger,B^\dagger)\\
p=\cd\ot b_j\ot\cd\ot b_1&\mapsto&p^\dagger=\cd\ot t(b_j)\ot\cd\ot t(b_1)
\end{eqnarray*}
{}From Theorem \ref{th:main} it suffices to show that 
this map is a bijection onto $\P^{(\la)}(\pb^\dagger,B^\dagger)$ such that 
$W(p)=\la+W(p^\dagger)$. Preservation of weight is immediate. To show
the bijectivity one has to notice that $\wt p^\dagger[j]-\wt p[j]$
does not depend on $j$. Thus one has $\wt p^\dagger[j]-\wt p[j]
=\wt p^\dagger-\wt p=-\la$, and hence
\[
\la_j(p^\dagger)=\la+\wt p^\dagger[j]=\wt p[j]=\vphi(b_{j+1})=\veps(b_j).
\]
Note that $p\in\PB_0$ (cf. Proposition \ref{prop:P_0} \& Corollary 
\ref{cor:useful}).
In view of (1) this equality concludes the bijectivity.
\qed

We now consider the $C^{(1)}_n$ case. For an odd positive integer $l$,
consider a finite crystal $B^{1,l}$ given by
\[
B^{1,l} = \left\{(x_1,\ldots, x_n,\overline{x}_n,\ldots,\overline{x}_1) 
\left| 
\begin{array}{l}
x_i, \overline{x}_i \in {\bf Z}_{\ge 0}\,\forall\, i=1,\cd,n \\
\sum_{i=1}^n(x_i + \overline{x}_i) \in \{l, l-2, \ldots, 1\}
\end{array}
\right.\right\}.
\]
The crystal structure of $B^{1,l}$ is given by
\begin{eqnarray*}
\et{0} b &=& \begin{cases}
(x_1-2,x_2,\ldots,\overline{x}_2,\overline{x}_1) 
&\mbox{ if }x_1 \ge \overline{x}_1+2, \\
(x_1-1,x_2,\ldots,\overline{x}_2,\overline{x}_1+1) 
&\mbox{ if }x_1 = \overline{x}_1+1,\\
(x_1,x_2,\ldots,\overline{x}_2,\overline{x}_1+2) 
&\mbox{ if } x_1 \le \overline{x}_1,
\end{cases}\\
\et{i} b &=& \begin{cases}
(x_1,\ldots,x_i+1,x_{i+1}-1,\ldots,\overline{x}_1) 
&\mbox{ if } x_{i+1} > \overline{x}_{i+1},\\
(x_1,\ldots,\overline{x}_{i+1}+1,\overline{x}_i-1,\ldots,\overline{x}_1)
&\mbox{ if } x_{i+1} \le \overline{x}_{i+1},
\end{cases}\\
\et{n} b &=& (x_1,\ldots,x_n+1,\overline{x}_n-1,\ldots,\overline{x}_1),\\
\ft{0} b &=& \begin{cases}
(x_1+2,x_2,\ldots,\overline{x}_2,\overline{x}_1) 
&\mbox{ if }x_1 \ge \overline{x}_1, \\
(x_1+1,x_2,\ldots,\overline{x}_2,\overline{x}_1-1) 
&\mbox{ if } x_1 = \overline{x}_1-1\\
(x_1,x_2,\ldots,\overline{x}_2,\overline{x}_1-2) 
&\mbox{ if } x_1 \le \overline{x}_1-2,
\end{cases}\\
\ft{i} b &=& \begin{cases}
(x_1,\ldots,x_i-1,x_{i+1}+1,\ldots,\overline{x}_1) &
\mbox{ if } x_{i+1} \ge \overline{x}_{i+1},\\
(x_1,\ldots,\overline{x}_{i+1}-1,\overline{x}_i+1,\ldots,\overline{x}_1)
&\mbox{ if } x_{i+1} < \overline{x}_{i+1},
\end{cases}\\
\ft{n} b &=& (x_1,\ldots,x_n-1,\overline{x}_n+1,\ldots,\overline{x}_1),\\
\end{eqnarray*}
where $b= (x_1, \ldots, x_n, \overline{x}_n,\ldots,\overline{x}_1)$
and $i = 1, \ldots, n-1$.
If some component becomes negative upon application, it should be understood 
as $0$. The values of $\veps_i,\vphi_i$ read
\[
\begin{array}{ll}
\veps_0(b)=\frac{l-s(b)}2+(x_1-\ol{x}_1)_+,\quad
&\vphi_0(b)=\frac{l-s(b)}2+(\ol{x}_1-x_1)_+,\\
\veps_i(b)=\ol{x}_i+(x_{i+1}-\ol{x}_{i+1})_+,
&\vphi_i(b)=x_i+(\ol{x}_{i+1}-x_{i+1})_+,\\
\veps_n(b)=\ol{x}_n,
&\vphi_n(b)=x_n.
\end{array}
\]
Here $s(b)=\sum_{i=1}^n(x_i+\ol{x}_i),(x)_+=\max(x,0)$ and $i=1,\cd,n-1$.
$B^{1,l}$ is a level $\frac{l+1}{2}$ non-perfect crystal.
Now for a fixed $l$ set $B=B^{1,l}$.
The minimal elements of $B$ are grouped  as 
$\Bmin = \bigsqcup_{\mu \in (P^+_{cl})_{\frac{l+1}{2}}} B^{\le \mu}$, 
where for $\mu = \mu_0 \La_0 + \cdots + \mu_n \La_n$.  The set  
$B^{\le \mu}$ is given by
\begin{eqnarray*}
B^{\le \mu} &=& 
 \{b^\mu_k \mid \mu_{k-1}>0, 1 \le k \le n \} \cup 
 \{b^\mu_{\overline{k}} \mid \mu_k > 0, 1 \le k \le n \},\\
b^\mu_k &=& (\mu_1,\ldots,\mu_{k-1}-1,\mu_k+1,\ldots,\mu_n,\mu_n,\ldots,
\mu_{k-1}-1,\ldots,\mu_1),\\
b^\mu_{\overline{k}} &=& 
(\mu_1,\ldots,\mu_k-1,\ldots,\mu_n,\mu_n,\ldots,\mu_1).
\end{eqnarray*}
Next consider $B^\dagger=B^{1,1}$ by taking $l$ to be $1$.
Setting 
$$
b^\dagger_k = (x_i = \delta_{i k}, \overline{x}_i = 0),\quad  
b^\dagger_{\overline{k}} = 
( x_i = 0, \overline{x}_i = \delta_{i k})
$$
for $1 \le k \le n$, one has
$$
(B^\dagger)^{\le \mu} = \{b^\dagger_k \mid \mu_{k-1} > 0, 1 \le k \le n \} \cup
\{b^\dagger_{\overline{k}} \mid \mu_k > 0, 1 \le k \le n\}
$$
for $\mu$ as above.
Define the map $t : \Bmin \rightarrow B^\dagger$ by 
$$t\vert_{B^{\le \mu}}: b^\mu_k \mapsto b^\dagger_k \qquad \mbox{ for }
k \in \{1, \ldots, n, \overline{n}, \ldots, \overline{1} \}.
$$

We are to show that this $t$ satisfies the conditions (1) -- (4) in 
Proposition \ref{prop:case-1}. For our purpose fix a dominant integral weight
$\la\in(\Pcl^+)_{\frac{l-1}2}$ and define $\pb=\cd\ot\bb_j\ot\cd\ot\bb_1$ by
\[
\bb_j=\begin{cases}
b^{\la+\La_i}_{\ol{i}}
&\mbox{ if }j\equiv i\,(\mbox{mod }2n)\mbox{ for some }i\,(1\le i\le n),\\
b^{\la+\La_{i-1}}_i
&\mbox{ if }j\equiv1-i\,(\mbox{mod }2n)\mbox{ for some }i\,(1\le i\le n).
\end{cases}
\]
Note that $\veps(b^{\la+\La_i}_{\ol{i}})=\vphi(b^{\la+\La_{i-1}}_i)=\la+\La_i,
\veps(b^{\la+\La_{i-1}}_i)=\vphi(b^{\la+\La_i}_{\ol{i}})=\la+\La_{i-1}$.
$\pb$ becomes a reference path. 
Let us check (1) -- (4) in Proposition \ref{prop:case-1}. (1),(2) and (4)
are straightforward. To check (3) one can use the formula for $H_{BB}$
in \cite{KKM} section 5.7. (In \cite{KKM} our non-perfect case is not
considered. However, the formula itself is valid. Since the formula in
\cite{KKM} contains some misprints, we rewrite it below.) 
\[
H_{B^{1,l}B^{1,l}}(b\ot b')
=\max_{1\le j\le n}(\theta_j(b\ot b'),\theta'_j(b\ot b'),
                    \eta_j(b\ot b'),\eta'_j(b\ot b')),
\]
\begin{eqnarray*}
\theta_j(b\ot b')&=&
\sum_{k=1}^{j-1}(\ol{x}_k-\ol{x}'_k)+\frac12(s(b')-s(b)),\\
\theta'_j(b\ot b')&=&
\sum_{k=1}^{j-1}(x'_k-x_k)+\frac12(s(b)-s(b')),\\
\eta_j(b\ot b')&=&
\sum_{k=1}^{j-1}(\ol{x}_k-\ol{x}'_k)+(\ol{x}_j-x_j)+\frac12(s(b')-s(b)),\\
\eta'_j(b\ot b')&=&
\sum_{k=1}^{j-1}(x'_k-x_k)+(x'_j-\ol{x}'_j)+\frac12(s(b)-s(b')),
\end{eqnarray*}
where $b=(x_1,\ldots,x_n,\overline{x}_n,\ldots,\overline{x}_1),
b'=(x'_1,\ldots,x'_n,\overline{x}'_n,\ldots,\overline{x}'_1)$.

Therefore, the isomorphism in Proposition \ref{prop:case-1} holds 
with notations above.

\subsection{Example 2} \label{subsec:ex-2}
We consider the $A_{n-1}^{(1)}$ case. Let $B^{1,l}$ be the crystal base of
the symmetric tensor representation of $U'_q(A_{n-1}^{(1)})$ of degree $l$. 
As a set it reads
\[
B^{1,l}=\{(a_0,a_1,\cd,a_{n-1})\mid a_i\in\Zn,\sum_{i=0}^{n-1}a_i=l\}.
\]
For convenience we extend the definition of $a_i$ to $i\in\Z$ by setting
$a_{i+n}=a_i$ and use a simpler notation $(a_i)$ for $(a_0,a_1,\cd,a_{n-1})$.
For instance, $(a_{i-1})$ means $(a_{n-1},a_0,\cd,a_{n-2})$.
The actions of $\et{r},\ft{r}$ ($r=0,\cd,n-1$) are given by 
\[
\et{r}(a_i)=(a_i-\delta_{i,r}^{(n)}+\delta_{i,r-1}^{(n)}),\quad
\ft{r}(a_i)=(a_i+\delta_{i,r}^{(n)}-\delta_{i,r-1}^{(n)}).
\]
Here $\delta_{ij}^{(n)}=1$ ($i\equiv j\mbox{ mod }n$), $=0$ (otherwise). 
If some component becomes negative upon application, it should be understood 
as $0$. The values of $\veps,\vphi$ read as follows.
\[
\veps((a_i))=\sum_{i=0}^{n-1}a_i\La_i,\quad
\vphi((a_i))=\sum_{i=0}^{n-1}a_{i-1}\La_i.
\]
Thus $\lev B^{1,l}=l$ and all elements are minimal.
We introduce a $\Z$-linear automorphism $\sigma$ on $\Pcl$ by
$\sigma\La_i=\La_{i-1}$ ($\La_{-1}=\La_{n-1}$). 

Now consider the finite crystal $B=B^{1,l}\ot B^{1,m}$ ($l\ge m$)
and set $B^\dagger=B^{1,m}$. From Lemma \ref{lem:tensor} (1) 
the level of $B$ is $l$. Fix two dominant integral weights
$\la=\sum_{i=0}^{n-1}\la_i\La_i\in(\Pcl^+)_{l-m},
\mu=\sum_{i=0}^{n-1}\mu_i\La_i\in(\Pcl^+)_m$. From $(\la,\mu)$ we define 
a path 
\[
\pb^{(\la,\mu)}(j)=(\la_{i+j}+\mu_{i+2j})\ot(\mu_{i+2j-1})\in B.
\]
{}From Lemma \ref{lem:tensor} (2) we see $\pb^{(\la,\mu)}(j)\in
\Bmin$ and by (\ref{eq:ot-eps}),(\ref{eq:ot-phi}) we obtain
$\veps(\pb^{(\la,\mu)}(j))=\sigma^j\la+\sigma^{2j}\mu
=\vphi(\pb^{(\la,\mu)}(j+1))$. 
Therefore $\pb^{(\la,\mu)}$ is a reference path.

We would like to show 
\begin{equation} \label{eq:ex2-toshow}
\P(\pb^{(\la,\mu)},B)\simeq B(\la)\ot\P(\pb^{(\mu)},B^\dagger)
\quad\mbox{as $P$-weighted crystals}
\end{equation}
with $\pb^{(\mu)}(j)=(\mu_{i+j})$. To do this, consider the following
map
\begin{eqnarray}
\P(\pb^{(\la,\mu)},B)_0&\longrightarrow&\P(\pb^{(\mu)},B^\dagger)
\label{eq:ex2-map}\\
p&\mapsto&p^\dagger\nonumber
\end{eqnarray}
given by $p^\dagger(j)=(b^{(j)}_{i-j+1})$ for 
$p(j)=(a^{(j)}_i)\ot(b^{(j)}_i)$. Note that $\pb^{(\la,\mu)}$ is sent
to $\pb^{(\mu)}$ under this map. By Theorem \ref{th:main} it suffices 
to check the following items:
\begin{itemize}
\item[(i)] The map (\ref{eq:ex2-map}) is a bijection onto 
           $\P^{(\la)}(\pb^{(\mu)},B^\dagger)$.
\item[(ii)] $\wt p-\wt p^\dagger=\la$.
\item[(iii)] $E(p)=E(p^\dagger)$.
\end{itemize}
Since $p\in\P(\pb^{(\la,\mu)},B)_0$, one obtains 
(cf. Lemma \ref{lem:tensor} (2), Proposition \ref{prop:P_0})
\begin{eqnarray}
&&\vphi_i((a^{(j)}_i))=a^{(j)}_{i-1}\ge b^{(j)}_i=\veps_i((b^{(j)}_i))
\label{eq:ex2-1}\\
&&\vphi_i(p(j))=a^{(j)}_{i-1}+b^{(j)}_{i-1}-b^{(j)}_i=a^{(j-1)}_i
=\veps_i(p(j-1))\label{eq:ex2-2}
\end{eqnarray}
for any $i,j$. Taking sufficiently large $J$ and using (\ref{eq:ex2-2}),
one has
\begin{eqnarray*}
\wt p^\dagger[j]&=&\sum_i b^{(J)}_{i-J+1}\La_i+\sum_{k=j+1}^J\sum_i
(b^{(k)}_{i-k}-b^{(k)}_{i-k+1})\La_i\\
&=&\sum_i(b^{(J)}_{i-J+1}-a^{(J)}_{i-J}+a^{(j)}_{i-j})\La_i\\
&=&\sum_i a^{(j)}_{i-j}\La_i-\la.
\end{eqnarray*}
Thus the condition $\veps(p^\dagger(j))\le\la_j(p^\dagger)$ is 
equivalent to saying $b^{(j)}_{i-j+1}\le a^{(j)}_{i-j}$ for any $i$,
which is guaranteed by (\ref{eq:ex2-1}). This proves (i). For (ii)
one only has to notice that $\wt p[j]=\vphi(p(j+1))=\sum_i a^{(j)}_i
\La_i$.

In order to prove (iii), we set
\begin{eqnarray*}
E_L^{diff}&=&\sum_{j=1}^L j\bigl\{
H_{BB}(((a^{(j+1)}_i)\ot(b^{(j+1)}_i))\ot((a^{(j)}_i)\ot(b^{(j)}_i)))\\
&&\hspace{12mm}
-H_{B^\dagger B^\dagger}((b^{(j+1)}_{i-(j+1)+1})\ot(b^{(j)}_{i-j+1}))\bigr\}.
\end{eqnarray*}
We can assume $(a^{(j)}_i)\ot(b^{(j)}_i)\in\Bmin$
for $1\le j\le L+1$. Under such assumption the isomorphism 
$B^{1,l}\ot B^{1,m}\simeq B^{1,m}\ot B^{1,l}$ sends $(a_i)\ot(b_i)$ to 
$(b_{i+1})\ot(a_i-b_{i+1}+b_i)$ \cite{NY}. Thus, from Proposition
\ref{prop:HofTensor} we have 
\[
H_{BB}(((a_i)\ot(b_i))\ot((a'_i)\ot(b'_i)))=b_0+a'_0+b'_0+
H_{B^\dagger B^\dagger}((b_i)\ot(b'_{i+1})).
\]
Let us recall the following formula for $H_{B^{1,m}B^{1,m}}$ 
(cf. \cite{KKM} section 5.1).
\[
H_{B^{1,m}B^{1,m}}((b_i)\ot(b'_i))=\max_{0\le j\le n-1}\bigl(
\sum_{k=0}^{j-1}(b'_k-b_k)+b'_j\bigr)
\]
{}From this one gets
\begin{eqnarray*}
&&H_{B^\dagger B^\dagger}((b^{(j+1)}_i)\ot(b^{(j)}_{i+1}))
-H_{B^\dagger B^\dagger}((b^{(j+1)}_{i-j})\ot(b^{(j)}_{i-j+1}))\\
&&\hspace{2.3cm}=\sum_{k=1}^j(b^{(j+1)}_{k-j-1}-b^{(j)}_{k-j}).
\end{eqnarray*}
Using above facts and (\ref{eq:ex2-2}) one obtains
\[
E_L^{diff}=\sum_{j=1}^L\sum_{k=0}^{j-1}a^{(L)}_{-k}
+L\sum_{k=0}^L b^{(L+1)}_{-k}.
\]
This completes (iii).
We have finished proving (\ref{eq:ex2-toshow}). It is also known
\cite{KMN2} that $\P(\pb^{(\mu)},B^{1,m})\simeq B(\mu)$.
Therefore we have 
\[
\P(\pb^{(\la,\mu)},B^{1,l}\ot B^{1,m})\simeq B(\la)\ot B(\mu)
\quad\mbox{as $P$-weighted crystals}.
\]

The multi-component version is straightforward. Consider the finite
crystal $B^{1,l_1}\ot\cd\ot B^{1,l_s}$ ($l_1\ge\cd\ge l_s\ge l_{s+1}=0$).
For $\la^{(i)}\in(\Pcl^+)_{l_i-l_{i+1}}$ ($1\le i\le s$) we 
define a reference path $\pb^{(\la_1,\cd,\la_s)}$ by 
\begin{eqnarray*}
&&\mbox{the $k$-th tensor component of }\pb^{(\la_1,\cd,\la_s)}(j)\\
&&\hspace{1cm}=(\la^{(k)}_{i+kj-k+1}+\la^{(k+1)}_{i+(k+1)j-k+1}
+\cd+\la^{(s)}_{i+sj-k+1}).
\end{eqnarray*}
Then we have
\[
\P(\pb^{(\la_1,\cd,\la_s)},B^{1,l_1}\ot\cd\ot B^{1,l_s})\simeq 
B(\la_1)\ot\cd\ot B(\la_s).
\]
The proof will be given elsewhere.

\vskip.5cm\noindent
{\bf Acknowledgements}
\vskip.3cm

The authors thank Professor Masaki Kashiwara and Professor Tetsuji Miwa
for useful discussions and kind interest.

\end{document}